\documentclass[12pt]{article}

\newif\iffrench\frenchfalse
\usepackage[utf8]{inputenc}
\usepackage[T1]{fontenc}

\usepackage[english]{babel}
\usepackage{lmodern}

\usepackage{amssymb}
\usepackage{amsthm}

\usepackage{bbold} 
\usepackage{dsfont}
\usepackage{mathtools}
\usepackage{thmtools}

\usepackage{varioref} 
\usepackage{hyperref} 
\hypersetup{
    colorlinks=true,
    linkcolor=blue,
    }

\usepackage{enumitem} 
\usepackage[babel=true]{csquotes} 
\usepackage{color}
\usepackage{setspace}
\usepackage{array}
\usepackage{layout}

\usepackage[margin=20pt,font=small,labelfont=bf]{caption} 

\usepackage[
  backend=biber,
  maxbibnames=99,
  citestyle=numeric-comp,
  giveninits=true,
  url=false,
  isbn=false,
  doi=false,
]{biblatex}
\newcommand{\smarge}[2]{\usepackage[top=#1,bottom=#1+1cm,left=#1-#2,right=#1]{geometry}}

\newtheorem{thm}{Theorem}[section]
\newtheorem{ppn}[thm]{Proposition}
\newtheorem{cor}[thm]{Corollary}
\newtheorem{lem}[thm]{Lemma}

\theoremstyle{remark}
\newtheorem{req}[thm]{Remark}

\newcommand{\ioo}[1]{\left(#1\right)}
\newcommand{\ifo}[1]{\left[#1\right)}
\newcommand{\iof}[1]{\left(#1\right]}

\newcommand{\pth}[1]{\left(#1\right)}
\newcommand{\cro}[1]{\left[#1\right]}
\newcommand{\acc}[1]{\left\{#1\right\}}
\newcommand{\abs}[1]{\left|#1\right|}

\newcommand{\floor}[1]{\left\lfloor#1\right\rfloor}

\newcommand{\with}{\hspace{1cm}\text{with}\hspace{1cm}}

\newcommand{\ie}{i.e. } 
\newcommand{\eg}{e.g. }

\newcommand{\esp}{\hspace{1cm}}
\newcommand{\comment}[1]{\hspace{0.5cm}\text{#1}\hspace{0cm}}

\newcommand{\tq}{\hspace{0.25cm}/ \hspace{0.25cm}}
\newcommand{\vg}{,\,}

\newcommand{\goq}{\geqslant}
\newcommand{\loq}{\leqslant}
\newcommand{\eps}{\varepsilon}
\newcommand{\ind}{\mathbb{1}}

\newcommand{\fig}[3]{\begin{figure}[ht]\begin{center}\includegraphics[width=#1cm]{#2}\end{center}\caption{#3}\end{figure}}

\newcommand{\de}{\,\mathrm{d}}
\newcommand{\dr}{\partial}

\newcommand{\Er}{\mathds{R}}
\newcommand{\Zed}{\mathds{Z}}

\newcommand{\prb}[1]{\mathds{P}\pth{#1}}
\newcommand{\Esp}[1]{\mathds{E}\cro{#1}}

\newcommand{\qt}[1]{#1\,,\hspace{1cm}}

\smarge{2.5cm}{0cm}

\renewbibmacro{in:}{}

\title{Large deviations and the emergence of a logarithmic delay in a nonlocal linearised Fisher-KPP equation}
\author{Nathanaël Boutillon\footnote{Email address: \url{nathanael.boutillon@inrae.fr}, Tel: +33 4 32 72 23 51}
\\ 
\footnotesize{BioSP, INRAE, 84914, Avignon, France}\\
\footnotesize{Aix Marseille Univ, CNRS, I2M, Marseille, France} 
}

\begin{document}

\maketitle
\begin{abstract}
  We study a variant of the Fisher-KPP equation with nonlocal dispersal. Using the theory of large deviations, we show the emergence of a \enquote{Bramson-like} logarithmic delay for the linearised equation with step-like initial data. We conclude that the logarithmic delay emerges also for the solutions of the nonlinear equation. Previous papers found very precise results for the nonlinear equation with strong assumptions on the decay of the kernel. Our results are less precise, but they are valid for all thin-tailed kernels. 
\end{abstract}


\section{Introduction}

We are concerned with the spreading properties of the solution of the Fisher-KPP equation with nonlocal dispersal
\begin{equation}\label{eq:nlin}
  \left\{
  \begin{aligned}
    \dr_tv(t,x)&=[J*v(t,x)-v(t,x)]+f(v(t,x))& \esp t>0\vg x\in\Er,\\
    v(0,x)&=\ind_{\iof{-\infty,0}}(x)&\esp x\in\Er,
  \end{aligned}
  \right.
\end{equation}
where \[J*v(t,x)=\int_{\Er}v(t,x-y)J(\text{d}y).\]
We shall look for the emergence of a \enquote{Bramson-like} logarithmic delay.
Here, the reaction term $f$ is monostable and satisfies the KPP condition, and the dispersal kernel $J$ is thin-tailed. These assumptions are precised and explained below. 

The Cauchy problem~\eqref{eq:nlin} is akin to the classical (local) Fisher-KPP equation~\cite{Fis37,KPP37}
\begin{equation}\label{eq:local}
  \left\{
  \begin{aligned}
    \dr_tw(t,x)&=\dr_{xx}w(t,x)+f(w(t,x))& \esp t>0\vg x\in\Er,\\
    w(0,x)&=\ind_{\iof{-\infty,0}}(x)&\esp x\in\Er,
  \end{aligned}
  \right.
\end{equation}
but the diffusion term $\dr_{xx}w$ is replaced by the nonlocal dispersal term $J*v-v$.  

In both Cauchy problems~\eqref{eq:nlin} and~\eqref{eq:local}, under suitable assumptions on the reaction term $f$, the maximum principle holds and implies that $0\loq v\loq 1$ and $0\loq w\loq 1$ (see~\cite{KPP37} for the local equation and~\cite{Yag09} for the nonlocal equation). Therefore, in biological models,  one can see $v(t,x)$ and $w(t,x)$ as the environment occupancy by an invading species at time $t$ and at location $x$. 
For example, the Cauchy problem~\eqref{eq:nlin} arises in~\cite{FouMel04} as the infinite particle limit of a stochastic process modelling the range expansion of a plant. In this model, which stems from~\cite{BolPac97,LawDie02}, the reaction term $f(v)$ corresponds to the demography (birth and death) of individuals subject to competition, while the dispersal term $J*v-v$ (which is the generator of a jump process) describes the dispersal of the seeds. As ${\frac{1}{2}\Delta}$ is the generator of Brownian motion, the local equation~\eqref{eq:local} usually models small and frequent movements, as in the model of Fisher~\cite{Fis37}. 
See also~\cite{Ske51,Tur98} for related models.

Although our primary interest lies in the Cauchy problem~\eqref{eq:nlin}, we shall focus on the \emph{linearised} equation
\begin{equation}\label{eq:lin}
  \left\{
  \begin{aligned}
    \dr_tu(t,x)&=[J*u(t,x)-u(t,x)]+ru(t,x)& \esp t>0\vg x\in\Er,\\
    u(0,x)&=\ind_{\iof{-\infty,0}}(x)&\esp x\in\Er,
  \end{aligned}
  \right.
\end{equation}
where $r:=f'(0)>0$. The behaviour of the linear equation seems simpler to understand. From the study of the solutions of~\eqref{eq:lin}, we shall deduce a similar (weaker) result on the solutions of~\eqref{eq:nlin}. 

Now, let us focus on known results about the local equation, which has already been widely studied. We say that the reaction term $f$ is \emph{monostable} when it satisfies
\begin{equation}\label{eq:monostable}
  f\in C^1([0,1]),\esp f'(0)>0,\esp f(0)=f(1)=0,\esp \forall u\in\ioo{0,1}\vg f(u)>0.
\end{equation}
We assume throughout that the reaction term $f$ is monostable.
The solution $w$ of the local Cauchy problem~\eqref{eq:local} is then known to propagate with a finite speed $c^*>0$ (\cite{AroWei75}, Theorems 4.1 and 4.3), in the sense that
\begin{equation}\label{eq:propagation}
  \begin{aligned}
  \lim_{t\to+\infty}w(t,c't)=0  &\comment{for all $c'>c^*$,}\\
  \lim_{t\to+\infty}w(t,c''t)=1 &\comment{for all $c''\in\ioo{0,c^*}$.}
  \end{aligned}
\end{equation}
The speed $c^*$ which satisfies~\eqref{eq:propagation} is called the \emph{critical speed}. A \emph{travelling wave} of speed $c$ is a solution $w$ of the local equation which can be written in the form ${w(t,x):=W(x-ct)}$, where $W:\Er\to[0,1]$ satisfies $W(-\infty)=1$ and $W(+\infty)=0$. The function $W$ is called the profile of the travelling wave, and has the shape of an interface between the \emph{invaded zone}, where $W\simeq 1$, and the \emph{non-invaded zone}, where $W\simeq 0$. 
The solution of the Cauchy problem~\eqref{eq:local} converges \enquote{in shape} to the unique positive travelling wave of critical speed $c^*$~\cite{KPP37,Uch78}. We thus say that an \emph{invasion front} appears in the solution.
    
Now, we might wonder where exactly the invasion front is: is it really \enquote{near} the travelling wave? More precisely, given any level $\rho\in\ioo{0,1}$, let $\theta_{loc}^{\rho}(t)$ be the largest position such that ${w(t,\theta_{loc}^{\rho}(t))=\rho}$. The goal is to compare $\theta_{loc}^{\rho}(t)$, which represents the position of the invasion front, to $c^*t$, which is the position of the travelling wave with minimal speed. 

The \emph{KPP condition} on the reaction term $f$ is defined by
\begin{equation}\label{eq:kpp}
  \qt{\forall u\in\cro{0,1}} f(u)\loq f'(0)u.
\end{equation}
The KPP condition means that the per-capita growth rate $f(u)/u$ is maximal at 0. In biological models, it implies that there is no Allee effect (e.g.,~\cite{Tur98}). When the KPP condition is satisfied, the critical speed is known to be $c^*=2\lambda$, with $\lambda=\sqrt{f'(0)}$~\cite{KPP37}.
Bramson, in~\cite{Bra78}, showed that if $f'(u)$ is maximal at $u=0$ (which implies that the {KPP condition} is satisfied), and if the initial condition is $\ind_{\iof{-\infty,0}}$, then the invasion front is located well behind the position $x=c^*t$. More precisely, he showed 
\begin{equation}\label{eq:shift_local}
  \theta^{\rho}_{loc}(t)=c^*t-s\ln(t)+O_{t\to+\infty}(1)\with s=\frac{3}{2\lambda},
\end{equation}
where $\lambda=\sqrt{f'(0)}$. The term $-s\ln(t)$ is called the (Bramson) logarithmic delay. The same year, Uchiyama~\cite{Uch78} proved a slightly less precise result, but for much more general initial conditions and more general reaction term: the KPP condition is almost sufficient. A few years later, keeping the assumption that $f'(u)$ is maximal at $u=0$, Bramson~\cite{Bra83} showed that in fact the term $O(1)$ can be turned into a term $C+o(1)$ for some constant $C\in\Er$ (which depends on the initial condition). Moreover, he extended again the class of initial conditions for which~\eqref{eq:shift_local} holds, and for which the solution converges to a travelling wave. His proofs use the so-called McKean representation with a branching Brownian motion~\cite{Sko64,McK75}. Lau~\cite{Lau85} gave another proof of the convergence of the solution to a travelling wave, with the same assumptions as~\cite{Bra83}, but using analytic techniques. Likewise, more recently, the emergence of the logarithmic delay in~\eqref{eq:shift_local} has been proved thanks to interesting analytic techniques~\cite{HamNol13}, for initial conditions with a support bounded from above and for reaction terms satisfying the KPP condition. For example, the techniques of~\cite{HamNol13} have been used in~\cite{NRR17,NRR19} to prove refinements of~\eqref{eq:shift_local} up to order $\frac{1}{\sqrt{t}}$ (which before had been found formally in~\cite{EvS00}).

Bramson and Uchiyama's results, with their extensions, induce that in the local equation, when the KPP condition is satisfied and when the initial condition decreases sufficiently fast, the shift between the position $\theta^{\rho}_{loc}$ of the invasion front and the position $c^*t$ is unbounded. Such a situation does not always occur. See the works of Rothe~\cite{Rot81} and of Fife and McLeod~\cite{FifMcL77} for the study of reaction terms which give rise to  bounded shifts. See also the work of Giletti~\cite{Gil22} for the study of general monostable reaction terms (that is, satisfying~\eqref{eq:monostable} but not necessarily~\eqref{eq:kpp}), which always give rise to unbounded shifts.

Much fewer articles have focused on nonlocal equations. The emergence of a logarithmic delay has been shown for other nonlocal variants of the classical Fisher-KPP equation. The case of a competition which is nonlocal in space has been studied in~\cite{Pen18} (using probabilistic arguments) and in~\cite{BHR20} (using analytic arguments). In the latter, the authors also showed that for a slowly-decaying competition kernel, the delay takes an algebraic form. The case of a phenotype-dependent equation with a competition which is nonlocal in phenotype and a diffusion rate which depends on the phenotype, the \enquote{cane-toad equation}, has been studied in~\cite{BHR17}. 

We now turn our attention to the results we have at hand for the {nonlocal} Cauchy problems of interest~\eqref{eq:nlin} and~\eqref{eq:lin}. Here, the nonlocality is on the movements, not on the competition. Throughout this work, we let $J$ be a measure on $\Er$ and we assume that $J$ satisfies the assumptions:
\begin{equation}\label{eq:cj}
  J(\ioo{0,+\infty})> 0,\esp J(\Er)=1,
\end{equation}
together with the assumption:
\begin{equation}\label{eq:cjthin}
  \text{There exists $L>0$ such that}\esp \int_{\Er}e^{L |x|}J(\mathrm{d} x)<+\infty.
\end{equation}
Assumption~\eqref{eq:cjthin} means that the dispersal kernel is \emph{thin-tailed}. It means that individuals usually do not move too far away from their origin, and it implies that travelling waves do exist. If the tail of $J$ is too heavy, the invasion is accelerating~\cite{KLvD96,Gar11}. When the invasion is accelerating, self-similar travelling waves cannot exist due to a flattening of the solutions~\cite{GHR17}. 

Under those assumptions~\eqref{eq:cj} and~\eqref{eq:cjthin} on the dispersal kernel $J$ and the monostability assumptions~\eqref{eq:monostable} on the reaction term $f$, the Cauchy problem is well-posed (see \eg the argument at the beginning of~\cite{Yag09}). If we assume further that $J$ is continuous, then there exists, as in the local case, a minimal speed for travelling waves, called the \emph{critical speed}~\cite{Sch80,CovDup07,CJS08}. For the nonlocal equation, we denote by $c$ the critical speed. When the KPP condition~\eqref{eq:kpp} is satisfied, the expression of~$c$ is explicit:
\begin{equation}\label{eq:def_c}
  c=\inf_{\lambda>0}\frac{M(\lambda)+r-1}{\lambda}\esp \with M(\lambda)=\Esp{e^{\lambda X}}. 
\end{equation}
By Proposition 1.1 in~\cite{Gra22}, under Assumptions~\eqref{eq:cj} and~\eqref{eq:cjthin}, there exists a unique $\lambda_r>0$ at which the infimum is reached:
\begin{equation}\label{eq:def_lamc}
  c=\frac{M(\lambda_r)+r-1}{\lambda_r}.
\end{equation}
Under mild additional assumptions on the kernel $J$, the critical speed $c$ satisfies the propagation property~\eqref{eq:propagation} for all nonzero solutions $v$ with an initial support bounded from above (replacing $w$ by $v$ in~\eqref{eq:propagation}). See~\cite{LPL05}, Theorem 3.2.

Recently, Graham~\cite{Gra22} showed that the logarithmic delay~\eqref{eq:shift_local} can also arise for the \emph{nonlocal} equation~\eqref{eq:nlin}. He assumes that $f$ satisfies the KPP condition. He also makes the following slightly technical assumption on the kernel $J$: there exist $M>0$ and $B>0$ such that for all~$x\goq 0$,
\begin{equation}\label{eq:graham_assumption}
J([x,x+M])\goq BJ(\ioo{x+M,+\infty}).
\end{equation}
This assumption holds if $J$ has a compact support, is exponential or Gaussian, but may fail for other thin-tailed kernels. For example, let
\[J:=\alpha\sum_{k=1}^{+\infty}e^{-\abs{k}}\pth{\delta_{k^2}+\delta_{-k^2}},\]
where $\alpha>0$ is a normalisation constant such that $J(\Er)=1$. Then, for all $M>0$, there exists  $x\goq 0$ such that $J([x,x+M])=0$ and $J(\ifo{x+M,+\infty})>0$; thus~\eqref{eq:graham_assumption} cannot be satisfied. Roquejoffre~\cite{Roq22}, assuming that $J$ has a compact support, was able to turn the term $O(1)$ into $C+o(1)$ (for some $C\in\Er$). The techniques in both works are close to those of~\cite{HamNol13} (Graham combines them with a probabilistic argument, and Roquejoffre uses a refinement as in~\cite{NRR17}). Our goal is to relax their conditions on $J$ and to treat the general case of thin-tailed kernels. A first step towards this goal is to work on the \emph{linear} Cauchy problem and, using the maximum principle, to find an upper bound for the position of the invasion front for the nonlinear Cauchy problem. 

One specificity of this work is the use of the large deviation theory. In the works cited above, there are essentially two kinds of proofs: those which use analytic tools, as in~\cite{Fis37,KPP37,FifMcL77,Rot81,Lau85,HamNol13,NRR19,Gil22,Roq22} and those which use probabilistic tools, as in~\cite{McK75,Uch78,Bra78,Bra83,Pen18}. Graham and Uchiyama combine both kinds of proofs.
The methods here are probabilistic. We shall use a Feynman-Kac representation of the solution. Then, we shall use the large deviation theory to estimate precisely the terms of the sum which will arise. In fact, it seems natural to use the large deviation theory in the problem we have raised, because we are mainly concerned with the extreme behaviour of individuals modelled by the solution of the Cauchy problem -- much as the large deviation theory is concerned with the extreme values of random processes.

Note however that the use of the large deviation theory is not new in the study of local or nonlocal Fisher-KPP equations. Freidlin~\cite{Fre85} uses a large deviation principle over the paths of a Brownian motion to determine the area covered by the solution in a heterogeneous environment; Evans and Souganidis~\cite{EvaSou89} used a Hamilton-Jacobi version of the large deviation principle to get the same results as Freidlin with analytic methods. A recent article~\cite{ADK23} deals with an equation close to ours, but with bistable or ignition reaction. The authors use the same Feynman-Kac representation as we do, the terms of which are also estimated thanks to the large deviation theory. With such estimates they are able to derive requirements on the initial condition so that the population persists. Their problem, therefore, is different from ours but their methods are very close.

Finally, it is worth to mention~\cite{AddRee09}, in which Addario-Berry and Reed deal with branching random walks in discrete time. In his work, Graham~\cite{Gra22} explains quickly how, from Addario-Berry and Reed's main result, one can deduce the emergence of the logarithmic delay~\eqref{eq:shift_local} for the solution of the nonlinear nonlocal Cauchy problem~\eqref{eq:nlin} for a \emph{restricted} class of monostable reaction terms and general dispersal kernels (he does not enter much into the details because he focuses on \emph{general} monostable reaction term). The connection between Addario-Berry and Reed's work and the Cauchy problem~\eqref{eq:nlin} follows from a McKean-like representation of the solution of~\eqref{eq:nlin} using a well-chosen branching random walk (instead of a branching Brownian motion). The jump law of the branching random walk is the kernel~$J$, and the branching law depends on the nonlinear reaction term. In the classical McKean representation, which holds for the local equation~\eqref{eq:local}, the branching random walk is replaced by a branching Brownian motion.

\section{Main results}

The first result is a proposition which gives a representation of the solution $u(t,x)$ of the linear Cauchy problem~\eqref{eq:lin}.
Let $J$ be a kernel which satisfies the hypotheses~\eqref{eq:cj} and~\eqref{eq:cjthin}. The kernel $J$ is therefore a probability law. We consider a sequence $(X_k)_{k\goq 1}$ of real independent and identically distributed random variables, following the law $J$. We define the random walk
\begin{equation*}
  S_n=\sum_{k=1}^nX_k.
\end{equation*}
The following proposition is a Feynman-Kac representation of the solution of the linear Cauchy problem~\eqref{eq:lin}. 
\begin{ppn}\label{ppn:rep}
  Let $(S_n)_{n\goq 1}$ be the random walk defined above.
  Let $u(t,x)$ be a solution of the linear Cauchy problem~\eqref{eq:lin}. Then, for all $t\in\ifo{0,+\infty}$, for all $x\in\Er$,
  \begin{equation}\label{eq:representation}
    u(t,x)=e^{(r-1)t}\sum_{n=0}^{+\infty}\frac{t^n}{n!}\prb{S_n\goq x}.
  \end{equation}
\end{ppn}

As we look closely at each term of the sum~\eqref{eq:representation}, we observe that when $x$ is near the position $ct$, there is a trade-off between the two factors $\frac{t^n}{n!}$ and $\prb{S_n\goq x}$:
\begin{itemize}
\item for $t\simeq n$, the factor $\frac{t^n}{n!}$ is large, but the probability  $\prb{S_n\goq x}$ is small;
\item for $t\gg n$, the probability $\prb{S_n\goq x}$ is large but the factor $\frac{t^n}{n!}$ is small;
\item for $t\ll n$, both the probability $\prb{S_n\goq x}$ and the factor $\frac{t^n}{n!}$ are small.
\end{itemize}
In fact, we shall understand in the following that for $x=ct+o(\sqrt{t})$, the dominant terms of the sum~\eqref{eq:representation} are located around a  position $n\simeq\alpha t$ proportional to $t$. Therefore, when we study $u(t,ct)$, we will have to deal with probabilities of the form $\prb{\frac{S_n}{n}\goq \frac{c}{\alpha}}$, where $\frac{c}{\alpha}$ is a constant. The theory of large deviation provides an interesting framework to estimate precisely those very small probabilities as $n\to+\infty$, see \eg the introduction of~\cite{DemZei10}.
These observations will allow us to prove the main result.
\begin{thm}\label{thm:main}
  Let $J$ be a kernel which satisfies the hypotheses~\eqref{eq:cj} and~\eqref{eq:cjthin}. Let $u$ be a solution of the linear Cauchy problem~\eqref{eq:lin}. Take $\rho\in\ioo{0,1}$ and denote by
  \[\sigma^{\rho}(t)=\sup\acc{x\in\Er\tq u(t,x)\goq \rho}\] the position of the level $\rho$ of $u$ at time $t$. We have
\begin{equation*}
  \sigma^{\rho}(t)=ct-s\ln(t)+O_{t\to+\infty}(1)\with s=\frac{1}{2\lambda_r},
\end{equation*}
where $\lambda_r$ is defined by~\eqref{eq:def_lamc}.
\end{thm}

\begin{req}
  Naturally, as we are currently considering the linear equation, we are not expecting to find exactly the same result as in Equation~\eqref{eq:shift_local}.
  However, we do get the same result with the local counterpart of the linear Cauchy problem~\eqref{eq:lin},
\begin{equation}\label{eq:local_lin}
  \dr_tz=\dr_{xx}z+rz,\esp r=f'(0),
\end{equation}
together with the initial condition $z(t,1)=\frac{1}{4\pi}e^{-x^2/4}$.
We then have a simple explicit expression for the solutions of Equation~\eqref{eq:local_lin},
$z(t,x)=\frac{e^{rx}}{\sqrt{4\pi t}}e^{-x^2/4t}$.
We now look for the position $\sigma^{\rho}_{loc}(t)$ such that ${z(t,\sigma^{\rho}_{loc}(t))=\rho}$. A short computation yields, as $t\to+\infty$,
\[\sigma^{\rho}_{loc}(t)=2\sqrt{r}t-\frac{1}{2\sqrt{r}}\ln(t)+O(1)=c^*t-\frac{1}{2\lambda}\ln(t)+O(1),\]
which implies the same result as in Theorem~\ref{thm:main} but for the local equation. Note also that the values $\lambda=\sqrt{f'(0)}$ and $\lambda_r$ play a symmetric role for, respectively, the local equation and the nonlocal equation.
\end{req}

Thanks to the maximum principle, we then deduce the following corollary about the nonlinear equation.

\begin{cor}\label{cor:main}
    Let $J$ be a kernel which satisfies the hypotheses~\eqref{eq:cj} and~\eqref{eq:cjthin}. Let $f$ be a reaction term which satisfies the monostability conditions~\eqref{eq:monostable} and the KPP condition~\eqref{eq:kpp}. Let $v$ be a solution of the nonlinear Cauchy problem~\eqref{eq:nlin}. Take $\rho\in\ioo{0,1}$ and denote by
  \[\theta^{\rho}(t)=\sup\acc{x\in\Er\tq v(t,x)\goq \rho}\] the position of the level $\rho$ of $v$ at time $t$. We have
\begin{equation*}
  \theta^{\rho}(t)\loq ct-s\ln(t)+O_{t\to+\infty}(1)\with s=\frac{1}{2\lambda_r},
\end{equation*}
where $\lambda_r$ is defined by~\eqref{eq:def_lamc}. 
\end{cor}

\begin{req}
This corollary is consistent with the main results of Graham and Roquejoffre~\cite{Gra22,Roq22}, which include the situation when $J$ has compact support (see above).  
\end{req}

Section~\ref{s:large_deviations} gives preliminary results about the theory of large deviations; we shall use these results in Subsection~\ref{ss:proof_main} to prove Theorem~\ref{thm:main}. Subsection~\ref{ss:proof_prop} is devoted to the proof of Proposition~\ref{ppn:rep}. Finally, Subsection~\ref{ss:proof_cor} is devoted to the proof of Corollary~\ref{cor:main}.

\section{\label{s:large_deviations}Preliminary results about large deviations}

We first focus on the useful notions from large deviations theory. Let $J$ be a dispersal kernel satisfying the hypotheses~\eqref{eq:cj} and~\eqref{eq:cjthin} and let $X$ be a real random variable following the law $J$. Let $D=\acc{\zeta\in\Er\tq \Esp{e^{\zeta X}}<+\infty}$. The \emph{cumulant generating function} of the law $J$ is defined at all $\zeta\in D$ by
\begin{equation*}
  \Lambda(\zeta)=\ln\pth{\Esp{e^{\zeta X}}}.
\end{equation*}
 As $J$ decreases at least exponentially fast at $\pm\infty$, the set $D$ contains a neighbourhood of 0. The Legendre-Fenchel transform of $\Lambda$ is denoted by $\Lambda^*$ and plays a fundamental role in the statement of the theorem of Bahadur-Rao. It is defined, for all $z\in\Er$, by
  \begin{equation*}
    \Lambda^*(z)=\sup_{\zeta\in D}\cro{\zeta z-\Lambda(\zeta)}\in\cro{0,+\infty}.
  \end{equation*}
  Figure~\ref{fig:legendre} gives a graphical interpretation of the Legendre-Fenchel transform of a (symmetric) convex function.

\fig{7}{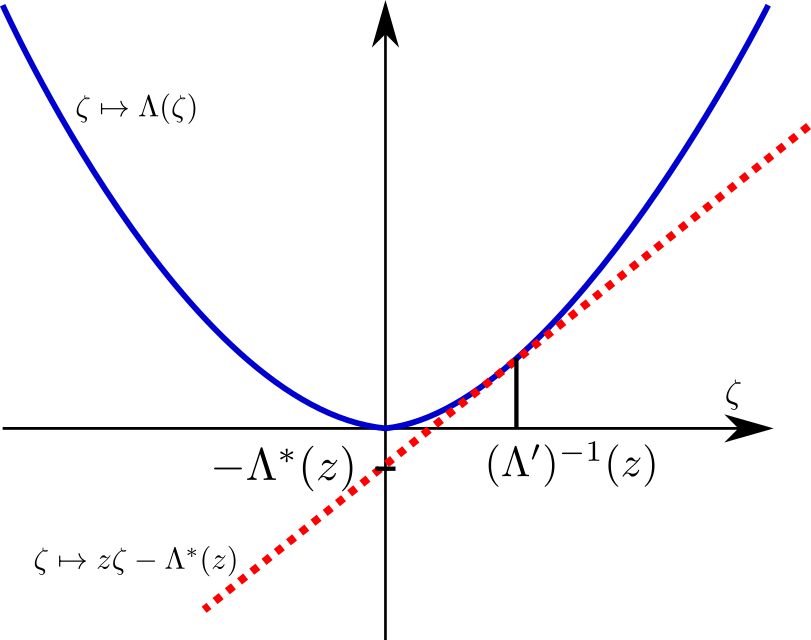}{\label{fig:legendre}Graphical interpretation of the Legendre-Fenchel transform of $\Lambda$ (inspired from~\cite{DemZei10}). The curve in blue looks like $\Lambda$, the dotted line in red is \emph{the} line tangent to the blue curve and with slope $z$. The dotted line in red crosses the vertical axis at $-\Lambda^*(z)$, the opposite of the Legendre-Fenchel transform of $\Lambda$ at $z$.}
  
  The following lemma contains general properties of the functions $\Lambda$ and $\Lambda^*$ (see~\cite{DemZei10}, Lemma 2.2.5).
  \begin{lem}\label{lem:prop_lambda}
    Assume that $J$ satisfies the hypotheses~\eqref{eq:cj} and~\eqref{eq:cjthin}. We then have:
    \begin{enumerate}
    \item The function $\Lambda$ is smooth (\ie $C^{\infty}$) over $D$. Let $A=\sup\acc{\text{\emph{supp}}(J)}$. The function $\Lambda^*$ is finite and smooth over $D^*=\ioo{-A,A}$;
    \item The functions $\Lambda$ and $\Lambda^*$ are nonnegative and strictly convex on, respectively, $D$ and $D^*$. Moreover, ${\Lambda(0)=\Lambda^*(0)=0}$;
    \item The function $\Lambda'$ is a bijection from $D$ to $D^*$. We define, for $z\in D^*$: ${\zeta_m(z):=(\Lambda')^{-1}(z)}$. We have, for all $z\in D^*$,
      \begin{equation}\label{eq:lamstar_explicite}
        \Lambda^*(z)=\zeta_m(z)z-\Lambda(\zeta_m(z))
      \end{equation}
      and
      \begin{equation*}
        (\Lambda^*)'(z)=\zeta_m(z).
      \end{equation*}
    \end{enumerate}
    
  \end{lem}
  

  \begin{proof}
    We sketch the proof of these classical properties (see also~\cite{DemZei10}, Lemma 2.2.5). 
    \begin{itemize}
    \item The nonnegativity of $\Lambda$ comes from Jensen's inequality. Since $J$ is 
      thin-tailed, $\Lambda$ is smooth in ${D}$;
    \item The strict convexity of $\Lambda$ follows: for $\zeta\in D$,
      \begin{align*}
        \Lambda''(\zeta)&=\frac{\Esp{X^2e^{\zeta X}}\Esp{e^{\zeta X}}-\Esp{Xe^{\zeta X}}^2}{\Esp{e^{\zeta X}}^2}\\
        &=\frac{\Esp{(Xe^{\zeta X/2})^2}\Esp{(e^{\zeta X/2})^2}-\Esp{(Xe^{\zeta X/2})(e^{\zeta X/2})}^2}{\Esp{e^{\zeta X}}^2}> 0,
      \end{align*}
      by the Cauchy-Schwarz inequality (the inequality is strict because $X$ is nonconstant);
    \item We have \[\Lambda'(\zeta)=\frac{\Esp{Xe^{\zeta X}}}{\Esp{e^{\zeta X}}}\underset{\zeta\to+\infty}{\to}\sup\acc{\text{supp}(J)},\] thus: ${\sup_{\zeta\in D}\Lambda'(\zeta)=\sup\acc{\text{supp}(J)}}$. Moreover, $\Lambda'$ is strictly increasing and smooth. Therefore, $\Lambda'$ is a smooth bijection from $D$ to $D^*$ with smooth converse $(\Lambda')^{-1}$;
    \item    For $z\in D^*$, the supremum in the definition of $\Lambda^*(z)$ is reached when $\Lambda'(\zeta)=z$. We get~\eqref{eq:lamstar_explicite}. We then deduce that $\Lambda^*$ is finite and smooth over $D^*$. Using~\eqref{eq:lamstar_explicite}, we get: $(\Lambda^*)'(z)=\zeta_m(z)$. Finally, $\zeta_m(z)$ is increasing, so $\Lambda^*$ is convex.
    \end{itemize}
  \end{proof}
  
  Let $(X_k)_{k\goq 1}$ be a sequence of independent and identically distributed random variables, following the law $J$, and define the corresponding random walk $(S_n)_{n\goq 1}$ as in the introduction.
  We are now ready to state, in our particular case, the theorem of Bahadur-Rao~\cite{BahRao60}. 
  This theorem will be useful in the proof of Theorem~\ref{thm:main}. Recall that, for $z\in D^*$, we have defined $\zeta_m(z)=(\Lambda')^{-1}(z)$. 

  \begin{thm}[Bahadur-Rao~\cite{BahRao60}]\label{thm:bahadur}
    Assume that $J$ is nonlattice. Take $M_1,M_2\in D^*$ with $M_2>M_1>0$. The random walk defined above satisfies
    \begin{equation}
      \prb{S_n\goq nz}=\frac{e^{-n\Lambda^*(z)}}{\zeta_m(z)\sqrt{2n\pi\Lambda''(\zeta_m(z))}}\pth{1+o(1)}\comment{as $n\to+\infty$,}
    \end{equation}
    uniformly in $z\in \cro{M_1,M_2}$.
  \end{thm}

The uniformity is not present in the original paper of Bahadur and Rao, but it has been proved by Petrov in~\cite{Pet65}. If $J$ is lattice (\ie $X_1$ is almost surely of the form $a+kb$, $k\in\Zed$, for some fixed $a,b\in\Er$), then the form is slightly different, but the proofs of the main results do not change.

\section{Proofs of the main results}

\subsection{\label{ss:proof_prop}Proof of Proposition~\ref{ppn:rep}}

\begin{proof}[Proof of Proposition~\ref{ppn:rep}]
If $u$ is a solution of the linear Cauchy problem~\eqref{eq:lin}, and if we set $\tilde{u}:=e^{-rt}u$, then $\tilde{u}$ is a solution of the linear equation ${\dr_t\tilde{u}=J*\tilde{u}-\tilde{u}}$ with the same initial condition.
Therefore, we may assume $r=0$.

Let $(Z_t)_{t\goq 0}$ be a Poisson process with rate $1$ and jump law $J$, and such that $Z_0=0$ almost surely. The infinitesimal generator of the process $(Z_t)_{t\goq 0}$ is \[\mathcal{G}f=J*f-f\] and, therefore, the function \[u(t,x):=\Esp{Z_t\goq x}\] solves $\dr_tu=J*u-u$ with initial condition $u(0,\cdot)=\ind_{\iof{-\infty,0}}$. Upon partitioning events according to the number of jumps made by the process in $[0,t]$, we conclude
\[u(t,x)=\sum_{n=0}^{+\infty}\frac{e^{-t}t^n}{n!}\prb{X_1+\hdots+X_n\goq x},\] where $\frac{e^{-t}t^n}{n!}$ is the probability that exactly $n$ jumps occurred in $[0,t]$. This yields the result for $r=0$. The conclusion follows by multiplying by $e^{rt}$.
\end{proof}

\begin{req}
    An important feature of Proposition~\ref{ppn:rep} is that it allows us to work on the discrete-time random walk $(X_1+\hdots+X_n)_{n\goq 0}=(S_n)_{n\goq 0}$ rather than on the continuous-time random walk $(Z_t)_{t\goq 0}$. In this situation, the theory of large deviations applies more easily. Another method to do the transformation from continuous time to discrete time is to consider the continuous-time random walk $(Z_t)_{t\goq 0}$ restricted to integer times, which gives the discrete-time random walk $(Z_n)_{n\goq 0}$. The jump law of the process $(Z_n)_{n\goq 0}$ is found by conditioning on the number of jumps made by the process $(Z_t)_{t\goq 0}$ between the times $0$ and $1$. Hence the sum arising in Proposition~\ref{ppn:rep} is directly incorporated into the jump law of $(Z_n)_{n\goq 0}$. This second method, which is inspired from the introduction of~\cite{AddRee09}, should lead to easier computations: one simply has to estimate the probabilities $\prb{Z_n\goq cn-\frac{1}{2\lambda_r}\ln(n)}$, and no more sum is involved. We shall rather concentrate on the first method, which seems more interesting from the point of view of modelling: indeed, it explicitly counts the jumps, \ie the generations. Such a record can be helpful if one wants to take into account the fact, for example, that mutations can arise at each generation.
\end{req}

\subsection{\label{ss:proof_main}Proof of Theorem~\ref{thm:main}}

Throughout Subsection~\ref{ss:proof_main}, we assume for convenience that $\sup\acc{\text{supp} J}=+\infty$, so that $\Lambda^*$ is defined on $D^*=\Er$ (the proof is almost the same if $\sup\acc{\text{supp} J}<+\infty$). We consider a nonnegative function $m(t)=o(\sqrt{t})$, and we note
\[x_t=ct-m(t).\] The function $m$ is intended to represent the delay, while $x_t$ is intended to represent the position of the front at time $t$. Our goal is to estimate $u(t,x_t)$ as $t$ grows to infinity. The idea of the proof is to cut the sum expressing $u(t,x_t)$ given by Proposition~\ref{ppn:rep},
\begin{equation}\label{eq:proof_sum_u}
  u(t,x_t)=e^{(r-1)t}\sum_{n=0}^{+\infty}\frac{t^n}{n!}\prb{S_n\goq x_t}
\end{equation}
into partial sums $S_{a,b}(t,x_t)$ of the form
\begin{equation*}
  S_{a,b}(t,x_t):=e^{(r-1)t}\sum_{n={\floor{at}+1}}^{\floor{bt}}\frac{t^n}{n!}\prb{S_n\goq x_t},
\end{equation*}
and to estimate each partial sum independently. We will show that there exists a value $\alpha>0$ such that the dominant terms of the total sum are located around $\floor{\alpha t}$. A large part of the proof is devoted to the estimation of those terms, that is, to the estimation of $S_{\alpha-\eps,\alpha+\eps}(t,x_t)$.

We begin with a lemma which estimates the partial sums $S_{a,b}(t,x_t)$. We use a version of the theorem of Bahadur-Rao (Theorem~\ref{thm:bahadur}) to turn the probability $\prb{S_n\goq x_t}$ into a more tractable expression. In the proofs, we will sometimes use the notation
\[p(t)=\Theta(q(t))\]
to say that there exist $0<K_-<K_+$ and $t_0> 0$ such that for all $t\goq t_0$, $K_-q(t)\loq p(t)\loq K_+q(t)$.

\begin{lem}\label{lem:est_psum}
  Let $b>a>0$. Define
  \[g({y}):={y}\pth{\ln\pth{\frac{e}{y}}-\Lambda^*\pth{\frac{c}{y}}}+r-1\]
  and
  \[h(t,y):=\frac{e^{(\Lambda^*)'(c/y)m(t)}}{t}e^{tg(y)}.\]
  There exist $K_-(a,b)\vg K_+(a,b)>0$ such that for $t$ large enough,
  \[K_-(a,b){\sum_{n={\floor{at}+1}}^{\floor{bt}}h\pth{t,\frac{n}{t}}}\loq S_{a,b}(t,x_t)\loq K_+(a,b){\sum_{n={\floor{at}+1}}^{\floor{bt}}h\pth{t,\frac{n}{t}}}.\]

\end{lem}

\begin{proof}
    We shall use the following property, resulting from Theorem~\ref{thm:bahadur}: for every subset $V\subset\ioo{0,+\infty}$, there exists $n_0(V)$ such that for all $n\goq n_0(V)$, for all $z\in V$,
\begin{align}\label{eq:property_bahadur}
  \frac{e^{-n\Lambda^*(z)}  }{2\zeta_m(z)\sqrt{2n \pi\Lambda''(\zeta_m(z))}}\loq\prb{S_n\goq n z}\loq 2\frac{e^{-n\Lambda^*(z)}  }{\zeta_m(z)\sqrt{2n \pi\Lambda''(\zeta_m(z))}}.
\end{align}
  Recall that $\zeta_m(z)$ is defined in Lemma~\ref{lem:prop_lambda}.

\paragraph{Step 1: Application of~\eqref{eq:property_bahadur}.}
Take $V=\cro{\frac{c}{2b},\frac{2c}{a}}$. Recall that $x_t=ct-m(t)$ where $m(t)=o(\sqrt{t})$ is nonnegative. Therefore, for $t$ large enough and for all integer $n$ between ${\floor{at}+1}$ and $\floor{b t}$, we have $x_t/n\in V$. Let $n_0(V)$ be the index given by Property~\eqref{eq:property_bahadur}. Consider $t_0\goq 0$ such that for all $t\goq t_0$, both ${\floor{at}+1}\goq n_0(V)$ and $x_t/n\in V$ hold. Take $t\goq t_0$ and apply Property~\eqref{eq:property_bahadur} to each integer $n$ between ${\floor{at}+1}$ and $\floor{b t}$, each time with $z_n=\frac{x_t}{n}\in V$. We find that there exist two constants $0<K_0<K_1$ independent of $t$ such that for each integer $n$ between ${\floor{at}+1}$ and $\floor{b t}$,
\[K_0\frac{e^{-n\Lambda^*\pth{\frac{x_t}{n}}}}{\zeta_m(x_t/n)\sqrt{2n\pi\Lambda''(\zeta_m(x_t/n))}}\loq\prb{S_n\goq x_t}\loq K_1\frac{e^{-n\Lambda^*\pth{\frac{x_t}{n}}}}{\zeta_m(x_t/n)\sqrt{2n\pi\Lambda''(\zeta_m(x_t/n))}}.\]
By Lemma~\ref{lem:prop_lambda}, the function $\zeta_m$ is continuous and positive over $\ioo{0,+\infty}$. Since $V$ is a compact subset of $\ioo{0,+\infty}$, we have \[0<\inf_{z\in V}\frac{1}{\zeta_m(z)\sqrt{2\pi\Lambda''(\zeta_m(z))}}\loq\sup_{z\in V}\frac{1}{\zeta_m(z)\sqrt{2\pi\Lambda''(\zeta_m(z))}}<+\infty.\]
Thus, summing over $n$, as $t\to+\infty$, 
\begin{equation}\label{eq:comparaison_ST}
   S_{a,b}(t,x_t)=\Theta\pth{T_{a,b}(t,x_t)}
\end{equation}
where
\[T_{a,b}(t,x_t):=e^{(r-1)t}\sum_{n=\floor{at}+1}^{\floor{bt}}\frac{t^n}{n!\sqrt{n}}e^{-n\Lambda^*\pth{\frac{x_t}{n}}}.\]
Our goal is now to estimate the sum $T_{a,b}(t,x_t)$.

\paragraph{Step 2: Estimation of $\Lambda^*$.} Take $t\goq t_0$ and an integer $n$ such that ${\floor{at}+1}\loq n\loq\floor{b t}$. Recall that $x_t=ct-m(t)$ where $m(t)=o(\sqrt{t})$. We have
  \begin{align*}
    \Lambda^*\pth{\frac{x_t}{n}}&=\Lambda^*\pth{\frac{ct}{n}-\frac{m(t)}{n}}=\Lambda^*\pth{\frac{ct}{n}}-(\Lambda^*)'\pth{\frac{ct}{n}}\frac{m(t)}{n}+O_{m/n\to 0}\pth{\frac{m(t)^2}{n^2}}.
  \end{align*}
  Thus there exists a time $t_1\goq t_0$ and a constant $K_2>0$ such that for all $t\goq t_1$ and for every integer $n$ such that ${\floor{a t}+1\loq n\loq\floor{b t}}$, we have
  \begin{align*}
    \abs{\Lambda^*\pth{\frac{x_t}{n}}-\pth{\Lambda^*\pth{\frac{ct}{n}}-(\Lambda^*)'\pth{\frac{ct}{n}}\frac{m(t)}{n}}}\loq K_2\frac{m(t)^2}{n^2}.
  \end{align*}
  We rewrite this as
  \begin{align}
    &\Lambda^*\pth{\frac{ct}{n}}-(\Lambda^*)'\pth{\frac{ct}{n}}\frac{m(t)}{n}-K_2\frac{m(t)^2}{n^2}\nonumber\\
    &\esp\loq\Lambda^*\pth{\frac{x_t}{n}}\nonumber\\
    &\esp\loq\Lambda^*\pth{\frac{ct}{n}}-(\Lambda^*)'\pth{\frac{ct}{n}}\frac{m(t)}{n}+K_2\frac{m(t)^2}{n^2}.\label{eq:estimate_step2}
  \end{align}
\paragraph{Step 3: Estimation of each single term of the sum $T_{a,b}(t,x_t)$.}  By Stirling's formula, as $n\to+\infty$, 
  \[n!=\Theta\pth{\sqrt{n}\exp(n\ln(n/e))}.\]
Using the second inequality in~\eqref{eq:estimate_step2}, we obtain a constant $K_3>0$ and a time $t_2\goq t_1$ such that, for all $t\goq t_2$, for all integer $n$ such that ${\floor{at}+1} \loq n\loq \floor{b t}$,
  \begin{align*}    
    &\frac{e^{(r-1)t}t^n}{n!\sqrt{n}}e^{-n\Lambda^*(x_t/n)}\\
    &\esp\goq\frac{K_3}{n}\exp\cro{n\pth{\ln\pth{e\frac{t}{n}}-\Lambda^*(ct/n)}+(r-1)t}\exp\cro{(\Lambda^*)'(ct/n)m(t)-K_2\frac{m(t)^2}{n}}\\
    &\esp=\frac{K_3}{n}\exp\cro{tg\pth{\frac{n}{t}}}\exp\cro{(\Lambda^*)'(ct/n)m(t)-K_2\frac{m(t)^2}{n}}.
  \end{align*}
  Since $m(t)=o(\sqrt{t})$, we have $m(t)^2=o(t)$. Therefore there exist $K_4>0$ and $t_3\goq t_2$ such that for all $t\goq t_3$, for all integer $n$ such that ${\floor{at}+1} \loq n\loq \floor{b t}$,
  \begin{align*}    
    \frac{e^{(r-1)t}t^n}{n!\sqrt{n}}e^{-n\Lambda^*(x_t/n)}
    &\goq\frac{K_4}{n}\exp\cro{tg\pth{\frac{n}{t}}}\exp\cro{(\Lambda^*)'(ct/n)m(t)}\goq\frac{K_4}{b}h\pth{t,\frac{n}{t}}.
  \end{align*}
  Thus, summing on all $n$ between ${\floor{at}+1}$ and $\floor{bt}$, we obtain an estimation of $T_{a,b}(t,x_t)$ which yields, together with Equation~\eqref{eq:comparaison_ST}, the existence of a constant $K_-(a,b)$ such that for all $t\goq t_3$,
  \[S_{a,b}(t,x_t)\goq K_-(a,b)\sum_{n={\floor{at}+1}}^{\floor{bt}}h\pth{t,\frac{n}{t}}.\]
  Using the first inequality in~\eqref{eq:estimate_step2}, we proceed to the same reasoning with reversed inequalities, and we obtain the existence of a constant $K_+(a,b)$ and a time $t'_3\goq t_2$ such that for all $t\goq t'_3$,
  \[S_{a,b}(t,x_t)\loq K_+(a,b)\sum_{n={\floor{at}+1}}^{\floor{bt}}h\pth{t,\frac{n}{t}}.\]
  The lemma is proved.
\end{proof}

Now, we introduce the value $\alpha$, which is constructed so that the most important terms of the sum~\eqref{eq:proof_sum_u} expressing $u(t,x_t)$ are located around the position $n\simeq \alpha t$.

\begin{lem}\label{lem:alpha}
   Set
  \begin{align*}
    \alpha=c\lambda_r-(r-1). 
  \end{align*}
  Then $\alpha>0$ and   \begin{equation}\label{eq:lemalpha1}
    \lambda_r=(\Lambda^*)'\pth{\frac{c}{\alpha}}.
  \end{equation}
 Moreover, the function $g$ defined in Lemma~\ref{lem:est_psum} is strictly concave on $\ioo{0,+\infty}$ and is maximal at $\alpha$. Finally, we have $g(\alpha)=0$, $g'(\alpha)=0$, and, for $y>0$, $y\neq\alpha$, we have $g(y)<0$. 
\end{lem}

\begin{req}
  We can write $h(t,y)$ in the form $h(t,y)=C(t,y)e^{tg(y)}$, where $C(t,y)$ is not too large as $t\to+\infty$. Hence, Lemmas~\ref{lem:est_psum} and~\ref{lem:alpha} entail that when $a<b<\alpha$ or $b>a>\alpha$, the partial sum $S_{a,b}(t,x_t)$ contains only terms which are exponentially small in $t$, thus the partial sum $S_{a,b}(t,x_t)$ is also exponentially small in $t$. Therefore, the dominant terms of the whole sum~\eqref{eq:proof_sum_u} expressing $u(t,x_t)$ are located around $n\simeq\alpha t$.
\end{req}

\begin{proof}[Proof of Lemma~\ref{lem:alpha}]
  Let $M(\lambda)=\Esp{e^{\lambda X}}$. Let $A(\lambda)=\frac{M(\lambda)+r-1}{\lambda}$. With this notation, ${c=\inf_{\lambda>0}A(\lambda)=A(\lambda_r)}$ and \[\alpha=c\lambda_r-(r-1)=M(\lambda_r).\] Hence $\alpha>0$. Since $\lambda_r$ is positive and $A$ is smooth and minimal at $\lambda_r$, we have ${A'(\lambda_r)=0}$, that is
  \begin{align*}
    \frac{1}{\lambda_r^2}\pth{\lambda_r M'(\lambda_r)-M(\lambda_r)-(r-1)}=0.
  \end{align*}
  Recall $\Lambda=\ln M$. Thus, differentiating,
  \[\lambda_r\Lambda'(\lambda_r)=\lambda_r\frac{M'(\lambda_r)}{M(\lambda_r)}=1+\frac{r-1}{M(\lambda_r)},\]
  which implies
  \[\Lambda'(\lambda_r)=\frac{1}{\lambda_r}+\frac{r-1}{\lambda_rM(\lambda_r)}=\frac{M(\lambda_r)+r-1}{\lambda_rM(\lambda_r)}=\frac{c}{M(\lambda_r)}=\frac{c}{\alpha}.\]
  Lemma~\ref{lem:prop_lambda} tells us that $(\Lambda^*)'=\zeta_m=(\Lambda')^{-1}$. Equality~\eqref{eq:lemalpha1} follows.

  We also have
  \begin{align*}
    \Lambda^*\pth{\frac{c}{\alpha}}&=\zeta_m\pth{\frac{c}{\alpha}}\frac{c}{\alpha}-\Lambda\pth{\zeta_m\pth{\frac{c}{\alpha}}}=\lambda_r\frac{c}{\alpha}-\Lambda(\lambda_r)
  \end{align*}
  so, using the fact that $\alpha=M(\lambda_r)$,
  \begin{equation}\label{eq:lemalpha2}
        \Lambda^*\pth{\frac{c}{\alpha}}+\ln(\alpha)-\frac{c\lambda_r}{\alpha}=0.
  \end{equation}
  With equalities~\eqref{eq:lemalpha1} and~\eqref{eq:lemalpha2} at hand, we are ready to conclude. Recall that \[g(y)=y\pth{\ln\pth{\frac{e}{y}}-\Lambda^*\pth{\frac{c}{y}}}+r-1.\]
  We have, by~\eqref{eq:lemalpha2},
  \[g(\alpha)=\alpha\pth{1-\ln(\alpha)-\Lambda^*(c/\alpha)}+r-1=\alpha\pth{1-\frac{c\lambda_r}{\alpha}}+r-1.\]
  Hence, since $\alpha=c\lambda_r-(r-1)$, we conclude that $g(\alpha)=0$. Furthermore, 
  \[g'(y)=-\ln(y)-\Lambda^*\pth{\frac{c}{y}}+\frac{c}{y}(\Lambda^*)'\pth{\frac{c}{y}}.\]
  Therefore, by~\eqref{eq:lemalpha1} and~\eqref{eq:lemalpha2}, we get $g'(\alpha)=0$. Finally, we have, for all $y>0$,
  \[g''(y)=-\frac{1}{y}-\frac{c^2}{y^3}(\Lambda^*)''\pth{\frac{c}{y}}<0.\]
  These elements allow us to conclude.
\end{proof}

\begin{lem}\label{lem:est_type_psum}
  The following asymptotic properties hold as $t\to+\infty$.
  \begin{enumerate}
  \item For $0<a<b<\alpha$, there exists $\chi>0$ such that
    \[S_{a,b}(t,x_t)=o(e^{-\chi t});\]
  \item For $B>A>\alpha$, there exists $\chi>0$ such that
    \[S_{A,B}(t,x_t)=o(e^{-\chi t});\]
  \item There exists $\eps>0$, two constants $K_-\vg K_+>0$ and a time $t_0>0$ such that for all $t\goq t_0$,
    \[K_-\frac{e^{\lambda_rm(t)}}{\sqrt{t}}\loq S_{\alpha-\eps,\alpha+\eps}(t,x_t)\loq K_+\frac{e^{\lambda_rm(t)}}{\sqrt{t}}.\] 
  \end{enumerate}
\end{lem}

\begin{proof}

Let $0<a<b<\alpha$ and let $\chi=-\frac{1}{2}\sup_{y\in\cro{a,b}}g(y)$. Then, by Lemma~\ref{lem:alpha}, $\chi$ is positive. From Lemma~\ref{lem:est_psum} and the fact that $m(t)=o(\sqrt{t})$, we deduce that the first point holds for this value of $\chi$. The same reasoning is valid for the second point as well. Now we prove the third point, which is more involved due to the equality $g(\alpha)=0$.
      
\paragraph{Step 1: Estimation of the partial sum $S_{\alpha-\eps,\alpha+\eps}(t,x_t)$; simple integrals arise.}
Take $\eps\in\ioo{0,\alpha/2}$. By Lemma~\ref{lem:est_psum}, there exist $K_0,K_1>0$ such that for $t$ large enough,
  \begin{align*}
     K_0t\pth{\frac{1}{t} \sum_{n=\floor{(\alpha-\eps)t}+1}^{\floor{(\alpha+\eps) t}}h\pth{t,\frac{n}{t}}}\loq S_{\alpha-\eps,\alpha+\eps}(t,x_t)\loq K_1t\pth{\frac{1}{t} \sum_{n=\floor{(\alpha-\eps)t}+1}^{\floor{(\alpha+\eps) t}}h\pth{t,\frac{n}{t}}}.
  \end{align*}
  By Lemma~\ref{lem:alpha}, there exists a neighbourhood $W$ of $\alpha$ such that for all $t>0$, the function ${y\mapsto h(t,y)}$ is decreasing on $W$. Up to reducing $\eps>0$, we may assume that $\cro{\alpha-2\eps,\alpha+2\eps}$ is included in $W$. Therefore, for all $t> 0$,
  \begin{align*}
    \frac{1}{t}\sum_{n=\floor{(\alpha-\eps)t}+1}^{\floor{(\alpha+\eps) t}}h\pth{t,\frac{n}{t}}\goq \sum_{n=\floor{(\alpha-\eps)t}+1}^{\floor{(\alpha+\eps) t}}\int_{\frac{n-1}{t}}^{\frac{n}{t}}h\pth{t,{y}}\de {y}
  \end{align*}
  and
  \begin{align*}
    \frac{1}{t}\sum_{n=\floor{(\alpha-\eps)t}+1}^{\floor{(\alpha+\eps) t}}h\pth{t,\frac{n}{t}}&\loq\sum_{n=\floor{(\alpha-\eps)t}+1}^{\floor{(\alpha+\eps) t}}\int_{\frac{n}{t}}^{\frac{n+1}{t}}h\pth{t,{y}}\de {y}.
  \end{align*}  
  Thus, there exist $K_2,K_3>0$ such that for $t$ large enough,
  \begin{equation}\label{eq:theta_sepsilon-}
    K_2{{t}\int_{-\eps}^{\eps}h(t,\alpha+y)\de y}\loq S_{\alpha-\eps,\alpha+\eps}(t,x_t)\loq K_3{{t}\int_{-\eps}^{\eps}h(t,\alpha+y)\de y}.
  \end{equation}

  \paragraph{Step 2: Estimation of the integral arising in~\eqref{eq:theta_sepsilon-} and conclusion.}   
  At the light of Equation~\eqref{eq:theta_sepsilon-}, we wish to estimate $h(t,\alpha+y)$ for $y$ close to $0$.
  By Lemma~\ref{lem:alpha}, we have $g''(\alpha)<0$ and
  \[g(\alpha+ y)=O(y^2).\]
  Therefore we can take $A',A''>0$, and reduce $\eps>0$ if necessary, so that for all $y\in\ioo{-\eps,\eps}$,
  \begin{equation*}
    -A'y^2\loq g(\alpha+ y)\loq -A''y^2.
  \end{equation*}
  Then, we have for some constants $K_5,K_6>0$ and for $t$ large enough,
  \begin{equation}\label{eq:theta_int_h}
    K_5\int_{-\eps}^{\eps}\frac{e^{(\Lambda^*)'\pth{\frac{c}{\alpha+y}}m(t)}}{t}e^{-tA'y^2}\de y\loq\int_{-\eps}^{\eps}h(t,y)\de y\loq K_6\int_{-\eps}^{\eps}\frac{e^{(\Lambda^*)'\pth{\frac{c}{\alpha+y}}m(t)}}{t}e^{-tA''y^2}\de y.    
  \end{equation}
  Now, note that \[(\Lambda^*)'\pth{\frac{c}{\alpha+y}}m(t)=\lambda_rm(t)+O_{y\to 0}(ym(t)).\]
  Thus, there are random variables $Z_t\sim \mathcal{N}\pth{O\pth{\frac{m(t)}{t}},\frac{1}{2tA'}}$ such that as $t$ goes to infinity,
  \begin{align*}
    \int_{-\eps}^{\eps}\frac{e^{(\Lambda^*)'\pth{\frac{c}{\alpha+y}}m(t)}}{t}e^{-tA'y^2}\de y&=\Theta\pth{\frac{e^{\lambda_rm(t)}}{t}\int_{-\eps}^{\eps}\exp\cro{-tA'{\pth{y+O\pth{\frac{m(t)}{t}}}^2}}\de y}\\
    &=\Theta\pth{\frac{{e^{\lambda_rm(t)}}}{t\sqrt{t}}\times\prb{-\eps<Z_t<\eps}}\\
    &=\Theta\pth{\frac{{e^{\lambda_rm(t)}}}{t\sqrt{t}}}.
  \end{align*}
  The last line holds because $m(t)=o(t)$, so that $Z_t$ converges to $0$ in probability, as $t\to+\infty$.
  Injecting twice this estimation into~\eqref{eq:theta_int_h} (once as such and once replacing $A'$ by $A''$), we deduce that there exist constants $K_7,K_8>0$ such that for $t$ large enough, 
  \begin{equation}\label{eq:minoration_int_h-}
    K_7\frac{{e^{\lambda_rm(t)}}}{t\sqrt{t}}\loq \int_{-\eps}^{\eps}h(t,y)\de y\loq K_8\frac{{e^{\lambda_rm(t)}}}{t\sqrt{t}}.
  \end{equation}
  Equation~\eqref{eq:minoration_int_h-}, together with Equation~\eqref{eq:theta_sepsilon-}, proves that the third point of the statement holds for the values of $\eps$ that we have selected in the beginning of Step 2.
\end{proof}

We are now ready to conclude the proof of the main theorem.

\begin{proof}[Proof of Theorem~\ref{thm:main}]
  Let $\eps>0$ be defined as in Lemma~\ref{lem:est_type_psum}. Take $a,B$ such that $0<a<\alpha<B$. By Proposition~\ref{ppn:rep}, we have
  \begin{align*}
    u(t,x_t)&=e^{(r-1)t}\sum_{n=0}^{+\infty}\frac{t^n}{n!}\prb{S_n\goq x_t},
  \end{align*}
  so, cutting the sum into five parts, we have
  \begin{align}
    u(t,x_t)&=\left(\sum_{n=0}^{\floor{at}}
    +\sum_{n=\floor{at}+1}^{\floor{(\alpha-\eps)t}}
    +\sum_{n=(\alpha-\eps)+1}^{\floor{(\alpha+\eps) t}}\right.\nonumber\\
    &\left.\esp\esp+\sum_{n=\floor{(\alpha+\eps)t}+1}^{\floor{Bt}}
    +\sum_{n=\floor{Bt}+1}^{+\infty}\right)\frac{e^{(r-1)t}t^n}{n!}\prb{S_n\goq x_t}.\label{eq:decomposition_u_parts}
  \end{align}
  Now, we estimate the first and last part of the sum (the other three will be estimated thanks to Lemma~\ref{lem:est_type_psum}).
  First, for $a>0$ small enough and $\chi_1>0$ small enough, we have, as $t\to+\infty$,
  \begin{align*}
    \sum_{n=0}^{\floor{at}}\frac{e^{(r-1)t}t^n}{n!}\prb{S_n\goq x_t}&\loq e^{rt}\prb{S_{\floor{at}}\goq x_t}=o(e^{-\chi_1 t}).
  \end{align*}
  The last estimation is obtained thanks to the Theorem~\ref{thm:bahadur} with $z=c/a$ (which becomes large when $a>0$ becomes small). 
  Second, we have, as $t\to+\infty$,
 \begin{align*}
   \sum_{n=\floor{Bt}+1}^{+\infty}\frac{e^{(r-1)t}t^n}{n!}\prb{S_n\goq x_t}&\loq\sum_{n=\floor{Bt}+1}^{+\infty}\frac{e^{(r-1)t}t^n}{n!}\\
   &\loq\frac{e^{(r-1)t}t^{\floor{Bt}}}{\floor{Bt}!}\sum_{n=0}^{+\infty}\frac{t^n}{n!}=e^{rt}\frac{t^{\floor{Bt}}}{\floor{Bt}!}=o(e^{-\chi_2 t})
 \end{align*}
 for $B>0$ large enough and $\chi_2>0$ small enough.

 Therefore, thanks to Lemma~\ref{lem:est_type_psum} and the decomposition~\eqref{eq:decomposition_u_parts}, we conclude that there exist constants $K_-,K_+>0$ and $\chi_3>0$ such that for $t$ large enough,
 \[K_-\frac{e^{\lambda_rm(t)}}{\sqrt{t}}\loq u(t,x_t)\loq K_+\frac{e^{\lambda_rm(t)}}{\sqrt{t}}+o(e^{-\chi_3 t}).\]
 Hence, there exists $t_0>0$ such that for all $t\goq t_0$,
 \[\frac{1}{2}K_-\frac{e^{\lambda_rm(t)}}{\sqrt{t}}\loq u(t,x_t)\loq 2K_+\frac{e^{\lambda_rm(t)}}{\sqrt{t}}.\]
 Hence, upon choosing $m(t)=\frac{1}{2\lambda_r}\ln(t)\pm C$ for a large $C$, we have:
 \begin{align*}
   u\pth{t,ct-\frac{1}{2\lambda_r}\ln(t)+C}&\ll\rho,\\
   u\pth{t,ct-\frac{1}{2\lambda_r}\ln(t)-C}&\gg\rho,
 \end{align*}
 as $t\to+\infty$. (Recall that $\rho$ is the level in which we are interested). Thus for $t$ large enough, $\abs{\sigma^{\rho}(t)-\frac{1}{2\lambda_r}\ln(t)}<C$. The conclusion of the theorem follows.
\end{proof}

\subsection{\label{ss:proof_cor}Proof of Corollary~\ref{cor:main}}

\begin{proof}[Proof of Corollary~\ref{cor:main}]
  We denote by $v$ the solution of the nonlinear Cauchy problem~\eqref{eq:nlin} and by $u$ the solution of the linear Cauchy problem~\eqref{eq:lin} with $r=f'(0)$. As the reaction term satisfies the KPP condition~\eqref{eq:kpp}, the function $v$ is a subsolution of the linear Cauchy problem~\eqref{eq:lin}. As $u$ and $v$ have the same initial condition, the maximum principle tells us that $v\loq u$. Finally, when we apply Theorem~\ref{thm:main} to the function $u$, we get the conclusion of Corollary~\ref{cor:main}.
\end{proof}

\section*{Acknowledgements}
This work was supported by the French Agence Nationale de la Recherche (ANR-18-CE45-0019 `RESISTE'), and by the Chaire Modélisation Mathématique et Biodiversité (École Polytechnique, Muséum national d’Histoire naturelle, Fondation de l’École Polytechnique, VEOLIA Environnement).

The author would like to express his gratitude to Raphaël Forien, François Hamel and Lionel Roques for their advice and constant support. He also would like to thank the Department of Mathematics at the Università di Roma \enquote{La Sapienza}, and the Department of Statistics at the University of Oxford, for their hospitality.

\printbibliography

@article{HamNol13,
	title = {A short proof of the logarithmic Bramson correction in Fisher-{KPP} equations},
	volume = {8},
	issn = {1556-181X},
	url = {http://aimsciences.org//article/doi/10.3934/nhm.2013.8.275},
	doi = {10.3934/nhm.2013.8.275},
	pages = {275--289},
	number = {1},
	journaltitle = {Netw. Heterog. Media},
	author = {Hamel, François and Nolen, James and Roquejoffre, Jean-Michel and Ryzhik, Lenya},
	urldate = {2021-04-15},
	date = {2013},
	langid = {english},
}

@article{Rot81,
  title={Convergence to pushed fronts},
  author={Rothe, Franz},
  journal={Rocky Mt. J. Math.},
  pages={617--633},
  year={1981},
  publisher={JSTOR}
}

@article{Lau85,
	title = {On the nonlinear diffusion equation of Kolmogorov, Petrovsky, and Piscounov},
	volume = {59},
	issn = {00220396},
	url = {https://linkinghub.elsevier.com/retrieve/pii/0022039685901378},
	doi = {10.1016/0022-0396(85)90137-8},
	pages = {44--70},
	number = {1},
	journaltitle = {J. Differ. Equations},
	author = {Lau, Ka-Sing},
	urldate = {2021-04-15},
	date = {1985},
	langid = {english},
}

@article{FifMcL77,
	title = {The approach of solutions of nonlinear diffusion equations to travelling front solutions},
	volume = {65},
	issn = {0003-9527, 1432-0673},
	url = {http://link.springer.com/10.1007/BF00250432},
	doi = {10.1007/BF00250432},
	pages = {335--361},
	number = {4},
	journaltitle = {Arch. Rational Mech. Anal.},
	author = {Fife, Paul C. and {McLeod}, J. B.},
	urldate = {2021-04-15},
	date = {1977},
	langid = {english},
}

@incollection{AroWei75,
	title = {Nonlinear diffusion in population genetics, combustion, and nerve pulse propagation},
	volume = {446},
	isbn = {978-3-540-07148-8},
	url = {http://link.springer.com/10.1007/BFb0070595},
	pages = {5--49},
	booktitle = {Partial Differential Equations and Related Topics},
	publisher = {Springer Berlin Heidelberg},
	author = {Aronson, D. G. and Weinberger, H. F.},
	date = {1975},
	langid = {english},
	note = {Series Title: Lecture Notes in Mathematics},
}

@article{Uch78,
	title = {The behavior of solutions of some non-linear diffusion equations for large time},
	volume = {18},
	issn = {2156-2261},
	url = {https://projecteuclid.org/journals/kyoto-journal-of-mathematics/volume-18/issue-3/The-behavior-of-solutions-of-some-non-linear-diffusion-equations/10.1215/kjm/1250522506.full},
	doi = {10.1215/kjm/1250522506},
	number = {3},
	journaltitle = {Kyoto J. Math.},
	author = {Uchiyama, Kōhei},
	date = {1978},
	pages = {453--508},
	langid = {english},
}

@article{KPP37,
	edition = {1},
	title = {Étude de l’équation de la diffusion avec croissance de la quantité de matière et son application à un problème biologique},
	pages = {1--26},
	journaltitle = {Bull. Univ. État Moscou},
	author = {Kolmogorov, Andrey N. and Petrovsky, Ivan G. and Piskunov, Nikolai S.},
	date = {1937},
	langid = {french},

}

@article{Fis37,
	title = {The wave of advance of advantageous genes},
	volume = {7},
	issn = {20501420},
	url = {http://doi.wiley.com/10.1111/j.1469-1809.1937.tb02153.x},
	doi = {10.1111/j.1469-1809.1937.tb02153.x},
	pages = {355--369},
	number = {4},
	journaltitle = {Ann. Eugenics},
	author = {Fisher, Ronald A.},
	date = {1937},
	langid = {english},
}

@article{Bra83,
	title = {Convergence of solutions of the Kolmogorov equation to travelling waves},
	volume = {44},
	url = {/paper/Convergence-of-solutions-of-the-Kolmogorov-equation-Bramson/c6177b6d18c334c32e38b19d404c3eeba022a9ba},
	journaltitle = {Mem. Am. Math. Soc.},
	author = {Bramson, Maury D.},
	urldate = {2021-04-15},
	date = {1983},
	langid = {english},
}

@article{Bra78,
	title = {Maximal displacement of branching Brownian motion},
	volume = {31},
	issn = {0010-3640},
	url = {https://onlinelibrary.wiley.com/doi/abs/10.1002/cpa.3160310502},
	doi = {10.1002/cpa.3160310502},
	pages = {531--581},
	number = {5},
	journaltitle = {Commun. Pure Appl. Math.},
	author = {Bramson, Maury D.},
	date = {1978},
}

@article{McK75,
	title = {Application of Brownian Motion to the Equation of Kolmogorov-Petrovskii-Piskunov},
	volume = {28},
	doi = {10.1002/cpa.3160280302},
	pages = {323--331},
	journaltitle = {Commun. Pure Appl. Math.},
	author = {{McKean}, Henry},
	date = {1975},
}

@article{CovDup07,
	title = {On a non-local equation arising in population dynamics},
	volume = {137},
	issn = {0308-2105, 1473-7124},
	url = {http://www.journals.cambridge.org/abstract_S0308210504000721},
	doi = {10.1017/S0308210504000721},
	pages = {727--755},
	number = {4},
	journaltitle = {P. Roy. Soc. Edinb. A},
	author = {Coville, Jerome and Dupaigne, Louis},
	date = {2007},
	langid = {english},
}

@article{EvS00,
	title = {Front propagation into unstable states: universal algebraic convergence towards uniformly translating pulled fronts},
	volume = {146},
	issn = {01672789},
	url = {https://linkinghub.elsevier.com/retrieve/pii/S0167278900000683},
	doi = {10.1016/S0167-2789(00)00068-3},
	shorttitle = {Front propagation into unstable states},
	pages = {1--99},
	number = {1},
	journaltitle = {Physica D},
	author = {Ebert, Ute and van Saarloos, Wim},
	urldate = {2021-04-15},
	date = {2000},
	langid = {english},
}

@book{DemZei10,
	location = {Berlin Heidelberg},
	edition = {2},
	title = {Large Deviations Techniques and Applications},
	isbn = {978-3-642-03310-0},
	url = {https://www.springer.com/gp/book/9783642033100},
	series = {Stochastic Modelling and Applied Probability},
	publisher = {Springer-Verlag},
	author = {Dembo, Amir and Zeitouni, Ofer},
	urldate = {2021-05-26},
	date = {2010},
	langid = {english},
	doi = {10.1007/978-3-642-03311-7},
}

@article{LPL05,
	title = {The effect of dispersal patterns on stream populations},
	volume = {47},
	pages = {749--772},
	number = {4},
	journaltitle = {{SIAM} Rev.},
	author = {Lutscher, Frithjof and Pachepsky, Elizaveta and Lewis, Mark A.},
	date = {2005},
}

@article{Pet65,
	title = {On the probabilities of large deviations for sums of independent random variables},
	volume = {10},
	pages = {287--298},
	number = {2},
	journaltitle = {Theor. Probab. Appl.},
	author = {Petrov, Valentin Vladimirovich},
	date = {1965},
	}

@article{Gra22,
	title = {The Bramson correction for integro-differential Fisher-{KPP} equations},
	volume = {20},
	issn = {1539-6746},
	url = {https://mathscinet.ams.org/mathscinet-getitem?mr=4374298},
	pages = {563--596},
	number = {2},
	journaltitle = {Commun. Math. Sci.},
	author = {Graham, Cole},
	date = {2022},
}

@article{BahRao60,
	title = {On deviations of the sample mean},
	volume = {31},
	pages = {1015--1027},
	number = {4},
	journaltitle = {Ann. Math. Stat.},
	author = {Bahadur, Raghu Raj and Rao, R. Ranga},
	date = {1960},
}

@article{KLvD96,
	title = {Dispersal Data and the Spread of Invading Organisms},
	volume = {77},
	rights = {© 1996 by the Ecological Society of America},
	issn = {1939-9170},
	url = {https://esajournals.onlinelibrary.wiley.com/doi/abs/10.2307/2265698},
	doi = {https://doi.org/10.2307/2265698},
	pages = {2027--2042},
	number = {7},
	journaltitle = {Ecology},
	author = {Kot, Mark and Lewis, Mark A. and Driessche, P. van den},
	urldate = {2021-05-09},
	date = {1996},
	langid = {english},
}

@article{AddRee09,
	title = {Minima in branching random walks},
	volume = {37},
	url = {https://hal.inria.fr/hal-00795281},
	pages = {1044--1079},
	journaltitle = {Ann. Probab.},
	author = {Addario-Berry, Louigi and Reed, Bruce},
	urldate = {2021-05-28},
	date = {2009},
	langid = {english},
}

@article{Gar11,
author = {Garnier, Jimmy},
title = {Accelerating Solutions in Integro-Differential Equations},
journal = {SIAM J. Math. Anal.},
volume = {43},
number = {4},
pages = {1955-1974},
year = {2011},
doi = {10.1137/10080693X},
}

@article{BolPac97,
	title = {Using moment equations to understand stochastically driven spatial pattern formation in ecological systems},
	volume = {52},
	pages = {179--197},
	number = {3},
	journaltitle = {Theor. Popul. Biol.},
	author = {Bolker, Benjamin and Pacala, Stephen W.},
	date = {1997},
}

@article{FouMel04,
	title = {A microscopic probabilistic description of a locally regulated population and macroscopic approximations},
	volume = {14},
	pages = {1880--1919},
	number = {4},
	journaltitle = {Ann. Appl. Probab.},
	author = {Fournier, Nicolas and Méléard, Sylvie},
	date = {2004},
}

@article{Fre85,
	title = {Limit theorems for large deviations and reaction-diffusion equations},
	pages = {639--675},
	journaltitle = {Ann. Probab.},
	author = {Freidlin, Mark},
	date = {1985},
}

@article{Ske51,
	title = {Random dispersal in theoretical populations},
	volume = {38},
	pages = {196--218},
	number = {1},
	journaltitle = {Biometrika},
	author = {Skellam, John Gordon},
	date = {1951},
}

@book{Tur98,
	title = {Quantitative analysis of movement: measuring and modeling population redistribution in animals and plants},
	shorttitle = {Quantitative analysis of movement},
	publisher = {Sinauer Associates},
	author = {Turchin, Peter},
	date = {1998},
}

@article {ADK23,
    AUTHOR = {Alfaro, Matthieu and Ducrot, Arnaud and Kang, Hao},
     TITLE = {Quantifying the threshold phenomenon for propagation in
              nonlocal diffusion equations},
   JOURNAL = {SIAM J. Math. Anal.},
  FJOURNAL = {SIAM Journal on Mathematical Analysis},
    VOLUME = {55},
      YEAR = {2023},
    NUMBER = {3},
     PAGES = {1596--1630},
      ISSN = {0036-1410,1095-7154},
   MRCLASS = {35B40 (35K57 45K05)},
  MRNUMBER = {4595400},
       DOI = {10.1137/22M1479099},
       URL = {https://doi.org/10.1137/22M1479099},
}

@article{LawDie02,
	title = {Moment approximations of individual-based models},
	pages = {252--262},
	journaltitle = {The Geometry of Ecological Interactions},
	author = {Law, Richard and Dieckmann, Ulf},
	date = {2002},
	note = {Cambridge Univ. Press},
}

@article{GHR17,
	title = {Transition fronts and stretching phenomena for a general class of reaction-dispersion equations},
	volume = {37},
	issn = {1078-0947},
	url = {https://mathscinet.ams.org/mathscinet-getitem?mr=3583498},
	doi = {10.3934/dcds.2017031},
	pages = {743--756},
	number = {2},
	journaltitle = {Discrete Contin. Dyn. Syst.},
	author = {Garnier, Jimmy and Hamel, François and Roques, Lionel},
	urldate = {2022-02-09},
	date = {2017},
	mrnumber = {3583498},
}

@article{Yag09,
	title = {Existence and nonexistence of traveling waves for a nonlocal monostable equation},
	volume = {45},
	issn = {0034-5318},
	url = {https://mathscinet.ams.org/mathscinet-getitem?mr=2597124},
	doi = {10.2977/prims/1260476648},
	pages = {925--953},
	number = {4},
	journaltitle = {Publ. Res. Inst. Math. Sci.},
	author = {Yagisita, Hiroki},
	urldate = {2022-02-09},
	date = {2009},
	mrnumber = {2597124},
}

@article{Sch80,
	title = {Travelling-front solutions for integro-differential equations. I},
	volume = {316},
	issn = {0075-4102},
	url = {https://mathscinet.ams.org/mathscinet-getitem?mr=581323},
	doi = {10.1515/crll.1980.316.54},
	pages = {54--70},
	journaltitle = {J. Reine Angew. Math.},
	author = {Schumacher, Konrad},
	urldate = {2022-02-09},
	date = {1980},
	mrnumber = {581323},
}

@article{NRR19,
	title = {Refined long-time asymptotics for Fisher–{KPP} fronts},
	volume = {21},
	pages = {1850072},
	number = {7},
	journaltitle = {Commun. Contemp. Math.},
	author = {Nolen, James and Roquejoffre, Jean-Michel and Ryzhik, Lenya},
	date = {2019},
}

@article{Pen18,
	title = {The spreading speed of solutions of the non-local Fisher–{KPP} equation},
	volume = {275},
	issn = {00221236},
	url = {https://linkinghub.elsevier.com/retrieve/pii/S0022123618303598},
	doi = {10.1016/j.jfa.2018.10.002},
	pages = {3259--3302},
	number = {12},
	journaltitle = {J. Funct. Anal.},
	author = {Penington, Sarah},
	date = {2018},
	langid = {english},
}

@article{BHR20,
	title = {The Bramson delay in the non-local Fisher-{KPP} equation},
	volume = {37},
	issn = {0294-1449},
	url = {https://www.sciencedirect.com/science/article/pii/S0294144919300757},
	doi = {10.1016/j.anihpc.2019.07.001},
	pages = {51--77},
	number = {1},
	journaltitle = {Ann. Inst. H. Poincaré, Anal. Non Linéaire},
	author = {Bouin, Emeric and Henderson, Christopher and Ryzhik, Lenya},
	date = {2020},
	langid = {english},
}

@article{BHR17,
	title = {The Bramson logarithmic delay in the cane toads equations},
	volume = {75},
	issn = {0033-569X},
	url = {https://mathscinet.ams.org/mathscinet-getitem?mr=3686514},
	doi = {10.1090/qam/1470},
	pages = {599--634},
	number = {4},
	journaltitle = {Quart. Appl. Math.},
	author = {Bouin, Emeric and Henderson, Christopher and Ryzhik, Lenya},
	date = {2017},
	mrnumber = {3686514},
}

@article{CJS08,
	title = {Nonlocal anisotropic dispersal with monostable nonlinearity},
	volume = {244},
	issn = {0022-0396},
	url = {https://mathscinet.ams.org/mathscinet-getitem?mr=2420515},
	doi = {10.1016/j.jde.2007.11.002},
	pages = {3080--3118},
	number = {12},
	journaltitle = {J. Differ. Equations},
	author = {Coville, Jérôme and Dávila, Juan and Martínez, Salomé},
	date = {2008},
	mrnumber = {2420515},
}

@article{NRR17,
	title = {Convergence to a single wave in the Fisher-{KPP} equation},
	volume = {38},
	pages = {629--646},
	number = {2},
	journaltitle = {Chin. Ann. Math. Ser. B},
	author = {Nolen, James and Roquejoffre, Jean-Michel and Ryzhik, Lenya},
	date = {2017},
}

@misc{Roq22,
	title = {Large time behaviour in nonlocal reaction-diffusion equations of the Fisher-{KPP} type},
	url = {http://arxiv.org/abs/2204.12246},
	journaltitle = {{arXiv}:2204.12246 [math]},
	author = {Roquejoffre, Jean-Michel},
	date = {2022},
	eprinttype = {arxiv},
	eprint = {2204.12246},
}

@article{Gil22,
	title = {Monostable pulled fronts and logarithmic drifts},
	volume = {29},
	issn = {1021-9722},
	url = {https://mathscinet.ams.org/mathscinet-getitem?mr=4412555},
	doi = {10.1007/s00030-022-00766-3},
	pages = {Paper No. 35},
	number = {4},
	journaltitle = {NoDEA Nonlinear Differential Equations Appl.},
	author = {Giletti, Thomas},
	date = {2022},
	mrnumber = {4412555},
}

@article{Sko64,
	title = {Branching Diffusion Processes},
	volume = {9},
	issn = {0040-585X, 1095-7219},
	url = {http://epubs.siam.org/doi/10.1137/1109059},
	doi = {10.1137/1109059},
	pages = {445--449},
	number = {3},
	journaltitle = {Theory Probab. Appl.},
	author = {Skorokhod, A. V.},
	date = {1964},
}

@article{EvaSou89,
	title = {A {PDE} approach to geometric optics for certain semilinear parabolic equations},
	volume = {38},
	number = {1},
	journal = {Indiana University mathematics journal},
	author = {Evans, Lawrence C. and Souganidis, Panagiotis E.},
	year = {1989},
	note = {Publisher: JSTOR},
}

\end{document}